\documentclass[a4paper]{scrartcl}
\usepackage[utf8]{inputenc}
\usepackage{amsfonts,amsthm,amssymb,amsmath, thmtools, thm-restate}
\usepackage[inline]{enumitem}
\usepackage{color}
\usepackage{graphicx} 
\usepackage{authblk}

\definecolor{mypink1}{RGB}{219, 48, 122}
\newtheorem{thm}{Theorem}
\newtheorem{lem}[thm]{Lemma}

\newtheorem{conj}[thm]{Conjecture}
\newtheorem{clm}{Claim}

\theoremstyle{remark}

\DeclareMathOperator{\Aut}{Aut}

\title{Distinguishing regular graphs from lists}
\author{Jakub Kwaśny, Marcin Stawiski\footnote{ Corresponding author; stawiski@agh.edu.pl}
}

\affil{AGH University,\\ Faculty of Applied Mathematics, \protect\\al. Mickiewicza 30, 30-059 Krakow, Poland}

\begin{document}
\maketitle
\begin{abstract}
An edge colouring of a graph is called \emph{distinguishing} if there is no non-trivial automorphism which preserves it. 
We prove that every at most countable, finite or infinite, connected regular graph of order at least $7$ admits a distinguishing edge colouring from any set of lists of length $2$. Furthermore, we show that the same holds for connected regular graphs of order $\kappa$ where $\kappa$ is a fixed point of the aleph hierarchy.
\end{abstract}

\section{Introduction}

All the colourings in this paper are (not necessarily proper) edge colourings, unless stated otherwise. We consider only simple graphs, which may be finite or infinite.
We say that $c$ \emph{breaks} an automorphism  of a graph  if there is an edge mapped to an edge of a different colour. A colouring  which breaks all the non-trivial automorphisms of a graph $G$ is called \emph{distinguishing}. The minimum number of colours in such a colouring is called the \emph{distinguishing index} of $G$, and it is denoted by $D'(G)$.

Distinguishing colourings were introduced by Babai \cite{BAB} in 1977 for vertex colourings, under the name \emph{asymmetric colourings}, which is still sometimes used. The  minimum number of colours in a distinguishing vertex colouring of a graph $G$ is called its \emph{distinguishing number} of $G$. The concept of breaking automorphisms played a crucial role in further research of Babai. They led to his proof of the existence of a quasi-polynomial algorithm for the graph isomorphism problem (see \cite{babaiisomorphism}). Distinguishing colourings gained quite a lot of attention since 1990s, and are still extensively studied. The most notable recent result in this area is the confirmation by Babai of the Infinite Motion Conjecture \cite{BABAI202264} proposed by Tucker \cite{tucker}. The edge variant of distinguishing colourings was introduced by Pilśniak and Kalinowski \cite{KP}.  Breaking automorphisms of graphs can be naturally generalised to breaking arbitrary group action on an arbitrary structure. These types of distinguishing colourings were studied first by Cameron, Neumann, and Saxl \cite{primitive} in 1984.

In this paper, we are interested in the list variant of edge distinguishing colourings of regular connected graphs. Erdős, Rubin, and Taylor \cite{erdos} started the study of list colourings in the setting of proper vertex colourings. The difference between the chromatic number of a graph $G$ and the required size of lists may be arbitrary large, even in the case where $\chi(G)=2$. In contrast, such examples are not known for  proper list edge colourings. The  List Colouring Conjecture \cite{bollobas-list-conjecture,haggkvist-list-conjecture} states that any graph $G$ has a proper edge colouring from any set of lists of length $\chi'(G)$.

Let $\mathcal{L}=\{L(e)\colon e\in E(G) \}$ be a set of lists for the edges of a graph $G$. The minimum length of each list in $\mathcal{L}$ such that there exists a distinguishing colouring from  $\mathcal{L}$   is called the \emph{list distinguishing index} of $G$, and it is denoted by $D'_l(G)$. Recently, Kwaśny and Stawiski \cite{KS_list_general} proved that each connected graph $G$ has the list distinguishing index at most $\Delta(G)-1$, with classified exceptions. 

\begin{thm}[Kwaśny, Stawiski \cite{KS_list_general}]\label{thm:delta-1}
 Let $G$ be a connected graph that is neither a double ray, a symmetric nor a bisymmetric
tree, $K_2$, $K_3$, $K_4$, $K_5$, $K_{3,3}$, $C_n$. Then
$D'_l(G) \le \Delta(G)-1$.
\end{thm}

There is a particular reason why distinguishing colourings of regular graphs may be of special interest. For any subgroup of the automorphism group of a given graph,  the orbit of a given vertex forms a transitive subgraph. This simple fact may be used in a construction of distinguishing colourings of arbitrary graphs because we can colour the connected components in the graph induced by the orbit of a given vertex with a small number of colours if this component is not isomorphic to $K_2$. This approach is not very useful in distinguishing vertex colourings because the distinguishing number of a regular connected graph may be arbitrarily large. Motivated by obtaining the bound for the distinguishing index of connected transitive graphs, Lehner, Pilśniak, and Stawiski \cite{LPSregular} studied the distinguishing colourings of regular connected graphs. They proved that each locally finite connected regular graph except $K_2$ has the distinguishing index at most $3$. This result was generalised by Grech and Kisielewicz to locally finite connected graphs for which $\delta(G)$ is at least half of $\Delta(G)$.  Moreover, Lehner, Pilśniak and Stawiski \cite{LPSregular} proposed the following conjecture.

\begin{conj}
\label{conj:mainregular}
If $G$ is a connected regular graph of order at least $7$, then $D'(G)\leq 2$.
\end{conj}

This was recently confirmed by Stawiski and Wilson \cite{Wilson} for connected regular graphs of infinite degree. We prove a stronger version of the Conjecture \ref{conj:mainregular} holds for locally finite connected graphs, namely the following theorem.

\begin{restatable}{thm}{main}
\label{thm:mainregular}
If $G$ is a connected locally finite regular graph of order at least $7$, then $D'_l(G)\leq 2$.
\end{restatable}

We also prove the same claim for a class of graphs which are not locally finite.
\begin{restatable}{thm}{mainxx}
\label{thm:aleph}
Let $\kappa$ be a cardinal number either equal to $\aleph_0$, or be a fixed point of the aleph hierarchy, i.e. $\kappa=\aleph_\kappa$. Let $G$ be a connected  $\kappa$-regular graph. Then $D'_l(G)\leq 2$.
\end{restatable}

In particular, it follows that the class of all cardinal numbers $\kappa$, for which every connected $\kappa$-regular graph has the list distinguishing index at most $2$, is unbounded, i.e. for every cardinal $\lambda$ there is a greater cardinal $\kappa$ with this property. In other words, the set of all such cardinals $\kappa$ forms a proper class. Moreover, for any regular cardinal $\kappa$, it is consistent with ZFC that every connected $2^{\kappa}$-regular graph has the list distinguishing index at most $2$. In particular, it holds for $\kappa=2^{\aleph_0}$. This follows from using the basic forcing method of Solovay, which is an extension of the classic method of Cohen, to force that $2^\kappa$ is a fixed point of the aleph hierarchy. We can obtain even stronger result by using the Easton forcing \cite{Easton}. This allows us to obtain a model of ZFC in which for every infinite regular cardinal $\kappa$, the cardinal $2^\kappa$ is a fixed point of the aleph hierarchy. Details of both methods may be found in a monograph of Jech \cite{Jech}.
\begin{thm}\label{thm:aleph_xD}
It is consistent with ZFC that for each infinite regular cardinal $\kappa$, every $2^{\kappa}$-regular connected graph has the list distinguishing index at most $2$.
\end{thm}
The relation between the distinguishing index and the list distinguishing index may also be of special interest.
It is clear that $D'(G)\leq D'_l(G)$. Examples of graphs for which $D'(G)<D'_l(G)$ are not known. Motivated by this, Kwaśny and Stawiski \cite{KS_list_general} proposed the following conjecture.
\begin{conj}\label{ourconjecture} Let $G$ be a connected graph without  a component isomorphic to $K_2$. Then $D'_l(G)=D'(G)$.
\end{conj}
The results in \cite{KS_list_general} imply that the statement of this conjecture holds for subcubic graphs. It follows from Theorems \ref{thm:mainregular} and \ref{thm:aleph} that this claim also applies to at most countable  connected regular graphs, and for $\kappa$-regular connected graphs where $\kappa$ is a fixed point of the aleph hierarchy.

The mentioned Infinite Motion Conjecture states that every locally finite graph such that every non-trivial automorphism moves infinitely many vertices has the distinguishing number equal to $2$. The edge variant of this conjecture was proposed by Broere and Pilśniak.
\begin{conj}[Broere, Pilśniak 2015 \cite{BP}]
    Let $G$ be a countable connected graph such that every non-trivial automorphism moves infinitely many edges. Then $D'(G)\leq 2$.
\end{conj}
This was confirmed by Lehner \cite{lehner-edge}. We propose the analogous conjecture for list colourings.
\begin{conj}
    Let $G$ be a countable connected graph such that every non-trivial automorphism moves infinitely many edges. Then $D'_l(G) \leq 2$.
\end{conj}
From the results of this paper, it follows that each countable regular graph satisfies the claim of this conjecture. It follows from the results of Kwaśny and Stawiski \cite{KS_list_general} that this conjecture is true in the case of subcubic graphs. The remaining cases are still open.

List distinguishing colourings in the setting of vertex colourings were introduced by Ferrara, Flesch and Gethner \cite{ferrara1} in 2011. They proposed a conjecture that the list distinguishing number is the same as the distinguishing number for every finite graph. There are only a few partial results towards this conjecture: for finite trees \cite{ferrara2}, finite interval graphs \cite{immel}, graphs with dihedral automorphism group \cite{ferrara1}, Cartesian product of two finite cliques \cite{furedi}, and Kneser graphs \cite{kneser}.
The known results about distinguishing graphs from lists by vertex colourings usually apply to very narrow classes of graphs, and they are usually less general than those obtained for non-list versions. In particular, none of the bounds for the list distinguishing number has been proved for infinite graphs. In contrast, this paper stands apart from this trend, as the obtained results are not only optimal, but also apply to general classes of graphs, including both finite and infinite graphs.

We would like to note that weaker results appeared on a paper by the same authors on arXiv \cite{KS_regular}. However, as we obtained stronger results before that paper was published, and we decided not to publish it outside of arXiv.

\section{Regular locally finite graphs}\label{main}

In this section, we prove the Theorem \ref{thm:mainregular}. We say that a colouring is \emph{near-distinguishing} if there is exactly one automorphism which preserves the colouring, and it interchanges some two vertices of a graph. In our proof, we use the following lemma. 

\begin{lem}
\label{lem:exceptions}
Let $G\in \{K_2, K_3, K_4, K_{3,3}, C_4, C_5\}$ and $\mathcal{L} = \{L_e\}_{e\in E(G)}$ be a set of lists of length $2$. Then, either $G$ admits a distinguishing edge colouring from these lists, or $G$ admits a near-distinguishing colouring.
\end{lem}
\begin{proof}
Assume that $G$ does not have a distinguishing edge colouring from the given set of lists. Note that $G$ is traceable. For each Hamiltonian path $P$, consider the value $s(P) = \sum_{e\in E(P)} \min(L_e)$ (without loss of generality, we can assume that the elements of the lists are natural numbers), and pick a path $P_0$ which minimises this sum. On that path, choose the lower values from the lists, and on all the other edges choose the higher values. Then, $P_0$ is the only Hamiltonian path with the sum of colours on the edges equal to $s(P_0)$, and therefore it is setwise fixed, and the only non-trivial automorphism that preserves the colouring is the one interchanging the end-vertices of $P_0$.
\end{proof}

Note that in the proof of Lemma \ref{lem:exceptions} we could maximise the sum of maxima, which in fact would give us a second, different near-distinguishing edge colouring from lists. We also need to know the list distinguishing index of finite complete graphs.

\begin{lem}
\label{lem:K_n}
If $n\ge 6$, then $D'_l(K_n) = 2$.
\end{lem}
\begin{proof}
We use the fact that $K_n$, for $n\ge 6$, contains an asymmetric spanning subgraph $H$. We proceed as in Lemma \ref{lem:exceptions}, considering all copies of $H$ and choosing one that minimises $s(H)$. Then again, on the edges of that copy we pick the lower values from the lists, and on all the other edges we pick the higher values from the lists. Consequently, there is exactly one copy $H_0$ of $H$ with the sum of colours on the edges equal to $s(H_0)$. Therefore, $H_0$ is setwise fixed, and since it is asymmetric, also pointwise fixed.
\end{proof}

We are ready to proceed with a proof of the main theorem.

\main*

Let $G=(V,E)$ be a regular graph of degree $\Delta \geq 2$ which is not an element of the exceptional family $\mathcal{H}=\{K_2, K_3, K_4, K_5, K_{3,3}, C_4, C_5\}$. The proof is by induction on $\Delta$. 

By Theorem \ref{thm:delta-1}  every connected 2-regular or 3-regular graph $G$ which is not in $\mathcal{H}$ has the list distinguishing index at most 2. This covers the base case of $\Delta\in \{2,3\}$. 

For the remaining part of the proof assume that $\Delta \geq 4$, and  that any $\Delta'$-regular graph with $\Delta'<\Delta$ has a list distinguishing edge colouring with $2$ colours unless it lies in $\mathcal{H}$. By Lemma \ref{lem:K_n}, we can also assume that $G$ is not a complete graph. Let $\mathcal{L} = \{L_e\}_{e\in E(G)}$ be any set of lists of length $2$, and let $C = \bigcup \mathcal{L}$ be the set of all available colours.

Consider a vertex $r\in V$. We will be using the automorphism group of a rooted graph $(G,r)$, which consists of exactly these automorphisms of $G$ which fix $r$. We denote this group by $\Aut(G,r)$. Our first task is to choose a starting vertex $r$ and colour the edges incident to $r$. 

Let $C=\{c_0, c_1, \dots\}$ be some arbitrary enumeration of the set of colours. For any vertex $v$, we count the number of occurrences of $c_0$ on the lists of the edges incident to $v$. Define $V_0 \subseteq V$ as the set of vertices with the maximum number $k_0$ of these occurrences. Now, for each vertex $v\in V_0$, we count the number of occurrences of $c_1$ on the list of these edges incident to $v$ which do not contain $c_0$. Define $V_1 \subseteq V_0$ as the set of these vertices in $V_0$ with the maximum number $k_1$ of these occurrences (among the vertices of $V_0$). We continue with this procedure and obtain a descending family of sets $V_0\supseteq V_1 \supseteq \ldots$ with a non-empty intersection $T=\bigcap_{i\ge 0} V_i$ (it is non-empty, as there are at most $\Delta$ indices $i$ such that $V_{i+1}\subsetneq   V_i$). 
We define $r$ as some vertex in $T$. 
Then, we colour $k_0$ edges incident to $r$ with the colour $c_0$, and $k_1$ edges incident to $r$ with the colour $c_1$, etc. Note that there is only one possible way to do it. We shall refer to the palette of $r$ we just constructed as the \emph{starting palette}. We shall ensure that this palette is never (or almost never) used again.

A frequent approach to the problem is constructing a colouring which breaks all the non-trivial automorphisms of $(G,r)$, and then arguing that $r$ must also be fixed by this colouring. Unfortunately, we cannot guarantee exactly that in our construction. However, we can construct a colouring that \emph{almost} fixes $r$, but instead satisfies an additional conditions. This additional conditions will allow us to perform one last modification of the colouring, which shall result in a distinguishing colouring of $G$. In the following theorem, we introduce the conditions, and then we prove that such modification is indeed possible.

\begin{thm} \label{clm: rooted} Let $r$ be a vertex chosen as above. There exists a list edge colouring $c$ of $G$, with the following properties: 
\begin{enumerate}
\item[(i)] $c$ breaks every non-trivial automorphism of $(G,r)$,
\item[(ii)] there exist exactly one or two vertices with the starting palette; one of them is the root $r$, and the other, if it exists, is a neighbour of $r$,
\item[(iii)] if there are two such vertices, then $T=V$, the palette of $r$ consists only of the colour $c_0$, and for each $k \in \{1, \ldots, \Delta \}$ the vertex $r$ has a neighbour with exactly $k$ incident edges with colour $c_0$. 
\end{enumerate}
\end{thm}

\emph{Proof of Theorem \textnormal{\ref{thm:mainregular}} using Theorem \textnormal{\ref{clm: rooted}}}.

Here we show how to construct a distinguishing colouring $c'$ of $G$ using the colouring $c$ from the claim. For brevity, we call $c_0$ red.

If $r$ is the only vertex with the starting palette, then it is fixed, therefore $c$ is a distinguishing colouring of $G$, and we put $c'=c$. Otherwise, there exists a vertex $x \in N(r)$ which has all the incident edges coloured red, and we can assume that there exists an automorphism moving $r$ to $x$. Let $y$ be the unique neighbour of $x$ with exactly one incident non-red edge.

\textbf{Case 1.} The vertex $y$ lies in $N(r)$. We argue that recolouring $yr$ with a colour other than red yields a distinguishing colouring $c'$ of $G$. First, notice that if we recolour any red edge incident to $r$, then the modified colouring is also distinguishing for $(G,r)$. Hence, after the recolouring, $x$ becomes fixed as the only vertex with all incident edges coloured red and $r$ is the only neighbour of $x$ with exactly one incident non-red edge.

\textbf{Case 2.} The vertex $y$ is not in $N(r)$. The vertex $y$ has a common non-red edge with exactly one vertex $z$. Let $k\ge 1$ be the number of blue edges incident to $z$. By the assumption that $\Delta \geq 4$, the vertex $r$ has at least two neighbours other than $x$ and $z$. Both of them have different numbers of incident non-red edges, so at least one of these numbers is different from $k-1$. Define $c'$ as obtained from $c$ by recolouring the edge between the corresponding vertex and $r$ with a colour other than red. Again, $x$ becomes the only vertex with all incident red edges, so it is fixed. By the previous argument, $c'$ is distinguishing for $(G,r)$, so we are left to confirm that $r$ is also fixed. However, $x$ has only two neighbours with exactly one incident non-red edge: $y$ and $r$, while $y$, unlike $r$, has a non-red edge to a vertex with exactly $k$ incident non-red edges. Therefore, $r$ cannot be mapped to $y$ and to any other vertex by an automorphism, so once more we obtain a distinguishing colouring of $G$. \hfill\qedsymbol

\vspace{3mm}

\emph{Proof of Theorem $\ref{clm: rooted}$.} 

We shall now construct a colouring that satisfies the assumptions of Theorem \ref{clm: rooted}. We shall describe an iterative procedure which, at each step, takes the \emph{first} not yet considered orbit with respect to some subgroup of $\Aut(G,r)$ and fixes (also with respect to $\Aut(G,r)$) each vertex in this orbit. This operation usually changes not only the considered orbit, but also some \emph{further} orbits, replacing them with their partitions. In order to maintain control over the entire algorithm and not leave any orbit, we need to specify carefully what the \emph{first} orbit is at any given point. This should be compatible with the partitioning of the orbits. Therefore, we now define a family of well-orderings of each possible set of orbits, i.e. for any possible subgroup of $\Aut(G,r)$. By well-ordering, we mean a linear ordering such that every non-empty subset has a least element. 

Let $\mathcal{S}_0$ be the set of orbits of vertices with respect to $\Aut(G,r)$. We establish some well-ordering $\leq_{\mathcal{S}_0}$  of $\mathcal{S}_0$ such that if $A,B \in \mathcal{S}_0$ and vertices of $A$ are closer to $r$ than $B$ (note that this distance is constant on any element of $\mathcal{S}_0$), then $A \leq_{\mathcal{S}_0} B$. Throughout the proof, we shall consider sets of orbits $\mathcal{S}$ with respect to nested subgroups of $\Aut(G,r)$. For each such set of orbits we establish a well-ordering $\leq_\mathcal{S}$ such that all these orderings are compatible in the following sense: if $\leq_{\mathcal{S}_1}$, $\leq_{\mathcal{S}_2}$ are the orderings of the sets of orbits corresponding to $\Gamma_1 \subseteq \Gamma_2 \subseteq \Aut(G,r)$, respectively, then for any $A_1,B_1 \in \mathcal{S}_1$, $A_1\subseteq A_2\in \mathcal{S}_2$, $B_1\subseteq B_2\in \mathcal{S}_2$, relation $A_2 \leq_{\mathcal{S}_2} B_2$ implies $A_1\leq_{\mathcal{S}_1} B_1$. Note that we shall choose such a family of well-orderings only for one nested chain of subgroups that will be iteratively constructed by our procedure. We shall drop the subscript in the ordering sign $\leq_{\mathcal{S}}$ whenever the related set $\mathcal{S}$ is clear.

Now, let $\mathcal{S}$ be one of such sets of orbits.
If we consider some $S \in \mathcal{S}$, then we refer to an edge between $S$ and vertex from $A\in \mathcal{S}$ as a \emph{back edge} if $A < S$, as a \emph{forward edge} if $S < A$, or as a \emph{horizontal edge} if $S=A$. 
Vertices connected to $S \in \mathcal{S}$ by forward or back edges will be referred to as \emph{forward neighbours} or \emph{back neighbours} of $S$, respectively.
Note that by definition of $\mathcal{S}$, every vertex in $S\in \mathcal{S}$ is incident to the same number $f$ of forward edges, the same number $b$ of back edges, and the same number $h$ of horizontal edges (and clearly $f+b+h=\Delta$). 
Further, note that by definition of $\mathcal{S}_0$ and the properties of $\leq_{\mathcal{S}_0}$ from the previous paragraph, $\bigcup \{S'\colon S' <S\}$ contains all the vertices that are closer to $r$ than $S$ (and perhaps some vertices at the same distance to $r$).
Thus, if $S\neq \{r\}$, then every vertex in $S$ has at least one back neighbour; i.e. $b \ge 1$.

We now proceed to the construction of the colouring $c$ which satisfies the conditions of Theorem \ref{clm: rooted}. 
Let $\mathcal{S}=\mathcal{S}_0$. It is clear that there can be only one orbit with $\Delta$ elements in $\mathcal{S}$ -- if there is such an orbit, then it must be the least element of $\mathcal{S} \setminus \{ \{ r\} \}$ and we denote it by $S_1$. If there is no orbit of size $\Delta$, then we do not define $S_1$. Note that if $S_1$ is defined, then the starting palette must consist only of the colour $c_0$.

We colour the edges of $G$ by iteratively running the procedure defined in the following paragraphs until all the vertices of $G$ are fixed with respect to the group $\Aut(G,r)$. The procedure takes the least orbit $S$ (with regard to the well-ordering $\le$) with at least one not coloured edge and colours all the edges incident to the vertices of $S$ which are not yet coloured, leaving $S$ fixed pointwise. The procedure needs some assumptions which we shall guarantee to remain satisfied after its completion.

\textbf{Assumptions for the procedure.} Let $S$ be the least element of $\mathcal{S}$ which has a vertex with an uncoloured incident edge ($S$ exists since $\le$ is a well-ordering). The following properties need to be satisfied before each execution of the procedure:
\begin{enumerate}[label = (A\arabic*)]
    \item every $S' <S$ has all incident edges coloured, no other edges are coloured, \label{itm:a1}
    \item the partial colouring defined so far fixes all $S'<S$ pointwise with respect to the group of automorphisms of $(G,r)$,\label{itm:a2}
    \item if $S'\in \mathcal{S}\setminus \{S_1\}$ and each element of $S'$ has exactly $b\geq 1$ coloured back edges, then $|S'|\leq \max \{1, \Delta-b \}$, \label{itm:gwiazdka}
    \item if a vertex $v$ has every incident edge coloured and its palette is the same as the starting palette, then either $v=r$ or $v$ is a unique neighbour of $r$ with this property. \label{itm:onlyonered}
\end{enumerate}
It is easy to verify that these assumptions are satisfied after the initial colouring of the edges incident to $r$.

\textbf{Procedure.}  
Each execution of the procedure consists of four phases. In phase one, we colour the horizontal edges if there are any. This gives us a distinguishing, or near-distinguishing, colouring of each component of $G[S]$ (the subgraph induced by $S$). In phase $2$, we choose a palette for each vertex in $S$. This choice will allow us to distinguish between the components of $G[S]$, and also complete the distinguishing of each component, if necessary. Phase $3$ involves colouring the forward edges of each vertex of $S$ with the previously selected palettes. At this stage, we partition the forward orbits of $S$ to guarantee the condition \ref{itm:gwiazdka}. Phase four is a slight modification of the colouring of the back edges to enforce the condition \ref{itm:onlyonered}.

\vspace{3mm}
\textit{Phase I: Horizontal edges}

Let $H$ be a component of the subgraph induced by the horizontal edges of $S$. Notice that $H$ is a connected $\Delta'$-regular graph for some $\Delta' < \Delta$. If $H \notin \mathcal{H}$, then we colour $H$ with a distinguishing colouring from the lists assigned to the edges of $H$, which exists by our induction assumption.
If $H \in \mathcal{H}$, then in each such component $H$, we choose a distinguishing or near-distinguishing edge colouring from the corresponding lists, given by Lemma \ref{lem:exceptions}. Additionally, if there is more than one component, and  the resulting colourings of all the components are isomorphic, then we recolour one of the components using the maxima in Lemma \ref{lem:exceptions} instead of the minima -- it will help us in the next phase. Denote by $P$ the Hamiltonian path used in Lemma \ref{lem:exceptions}. 

If $f=0$, then $h>0$ would imply (using (A3)) that $|H|>h=\Delta-b \ge |S|\ge |H|$, a contradiction. Therefore, $h=0$ and, again by (A3), $|S|=1$, so we are done. Note that it is not possible that $f=0$ for $S=S_1$ since $G$ is not a complete graph.

\vspace{3mm}
\textit{Phase II: Choosing forward palettes}

If $f>0$, then we now choose the forward palettes for the vertices in $S$. Let $A$ be the set of automorphisms of $\Aut(G,r)$ which preserve the current colouring (i.e. the one after completing Phase I).


For each vertex $x$ in $S$ we choose a number $i(x)$ between $0$ and $f$. This number will determine, how many orbits the vertex $x$ will split. Additionally, in the next phase, we shall guarantee that the vertices with different values of $i(x)$ shall receive distinct palettes. We will try to keep the values of $i(x)$ as close to $\lceil f/2 \rceil$ as possible, but sometimes we shall assign different values in order to ensure some additional properties.

In the following claim, the first two conditions will oblige us to break the remaining automorphisms interchanging the vertices of $S$. The next two conditions will keep us apart from breaking \ref{itm:onlyonered}. The last three conditions guarantee that the values of $i(x)$ are close enough to $\lceil f/2 \rceil$ to split sufficiently many orbits and satisfy \ref{itm:gwiazdka}. 

\begin{clm} \label{clm: i(x)} Let $H$ be a graph isomorphic to every component of $G[S]$. Assume that $f>0$. There exists an assignment $i\colon  S\rightarrow \{0,1,\ldots,f\}$ satisfying the following properties:
\begin{enumerate}[label = \textnormal{(B\textnormal{\arabic*})}]
    \item if $g\in A$ is an isomorphism which maps a component $H'$ of $G[S]$ into a different component of $G[S]$, then $i(g(s)) \neq i(s)$ for some $s \in H'$,
    \item if $H \in \mathcal{H}$ and both $s_1$ and $s_2$ are distinct end-vertices of the Hamiltonian path $P$ in a component of $G[S]$, then $i(s_1) \neq i(s_2)$,
    \item if $f\ge 2$, then only one vertex $s \in S$ may have $i(s)=0$ and it may happen only if $|S|\geq \Delta-1$, or $|S|=1$, or $f=2$ and $H \in \mathcal{H}\cup \{K_1\}$,
    \item if $f\ge 2$, then only one vertex $t \in S$ may have $i(t)=f$ and it may happen only if $|S|=\Delta$ (so $S=S_1$) and either $H=K_1$, or $H=K_2$ and $\Delta\in \{4,6\}$, \label{itm:b4}
    \item if $f\ge 2$ and there exists $x,y \in S$ such that $x$ and $y$ do not have the same sets of forward neighbours,  and $i(s)=0$, $i(t)=f$ for some $s,t \in S$, then $s$ and $t$ do not have the same sets of forward neighbours, \label{itm:b5}
    \item if $f\ge 5$, then for each $k=1,2,\ldots,\lceil\frac{f}2\rceil-2$, there are at most $k+1$ vertices with $i(s)\le k$ and at most $k+1$ vertices with $i(s)\ge f-k$ (and at most $k$ vertices with $i(x)\ge f-k$ when $S\neq S_1$). Moreover, if $|S|<f-2k$, then there are no vertices with $i(s)\le k$ and no vertices with $i(s)\ge f-k$.
    \item there exists a vertex $s\in S$ with $i(s)=\lceil \frac{f}2 \rceil$.
\end{enumerate}
\end{clm}
\begin{proof} 

Note first that it follows from \ref{itm:gwiazdka} and the fact that $|H| \geq h+1$ that the number of components in $G[S]$ is at most $\frac{\Delta-b}{h+1} = \frac{f+h}{h+1} = 1+\frac{f-1}{h+1} \le f$ if $S \neq S_1$ or at most $\frac{\Delta}{h+1} = \frac{f+h+1}{h+1} = 1+\frac{f}{h+1} \le f+1$ if $S=S_1$. 
Furthermore, if there are $f+1$ components of $G[S]$, then $h=0$ and $f=\Delta-1$, and if there are $f\geq 2$ components of $G[S]$, then either $h=0$, or $h=1$ and $f=2$. 

\textbf{Case 1.} $f=1$. Note that there is exactly one component $H$ of $S$, as otherwise $\Delta=f+1=2$. For $H \in \mathcal{H}$, we assign $i(s)=1$ to one arbitrarily chosen end-vertex $s$ of the blue Hamiltonian path $P$ and 0 to each remaining vertex of $H$. For $H \notin \mathcal{H}$, we just assign $i(s)=1$ to each vertex $s$ of $H$. 

\textbf{Case 2.} $f = 2$ and there is exactly one component $H$ of $S$. If $H \in \mathcal{H}$, then we put $i(s)=1$ everywhere except one end-vertex of the Hamiltonian path $P$, where we put $i(s)=0$. Otherwise, we just assign $1$ to each vertex of $H$.

\textbf{Case 3.} $f = 2$ and there are at least two components of $S$. Then, there are exactly two components, as otherwise $\Delta=f+1=3$. By the previous considerations about the number of components, we have either $h=0$, which means that $G[S] = 2K_1$, or $h=1$, which implies that $\Delta=4$, $S=S_1$, and $G[S]=2K_2$. The former is an exception in (B3) and we assign 1 to one vertex and 0 to the other. The latter is an exception in both (B3) and (B4) and we are here allowed to use $i(x)=0$ and $i(x)=2$, each only once. At the same time, we have enough freedom to satisfy condition (B5). 

\textbf{Case 4.} $f\ge 3$. If each component is a $K_1$, then we take any injective mapping into $\{0,1,\ldots,f-1\}$ (if $S\neq S_1$) or into $\{0,1,\ldots,f\}$ (if $S=S_1$), keeping the values as close to $\frac{f}2$ as possible, and satisfying (B5). Otherwise, consider any component $H$. Let $x$ and $y$ be the end-vertices of the previously chosen Hamiltonian path $P$, if $H \in \mathcal{H}$, or choose $x$ and $ y$ arbitrarily, otherwise. Consider all unordered pairs of different indices in $\{1,2,\dots,f-1\}$. There are $\binom{f-1}2$ such pairs which is at least $f-1$ whenever $f\geq 4$. On the other hand, since $h>0$ and $f>2$, the number of components is at most $f$. One of these components $H_0$ has a different colour of the Hamiltonian path $P$, so we may repeat a pair there. Therefore, we can use a unique pair on the vertices $x,y$ for each component $H\neq H_0$, starting with the pairs closest to $i_0=\lceil\frac{f}2\rceil$ (i.e. use an index only if all the pairs with both indices closer to $i_0$ have already been used), and repeat the pair $i_0, i_0-1$ on $H_0$. The remaining vertices receive $i_0$. 

If $f=3$ and $S\neq S_1$, then there are at most two components of $G[S]$ and only one pair of different indices in $\{1,2\}$. There is no problem if there is only one component. If $H \in \mathcal{H}$, then the number of pairs is also sufficient, because of the colours of $P$. Otherwise, each component has at least three vertices, and we simply change the index of a vertex other than $x$, $y$ from $i_0=2$ to 1.

If $f=3$ and $S=S_1$, then there are at most three components of $G[S]$ and still only one pair of different indices in $\{1,2\}$. If there are at most two such components, then we proceed as in the previous paragraph, therefore assume that there are three components. If $|H|\geq 4$, then we can again change the index of one vertex from 2 to 1 in one component, and repeat this with two vertices in another component. If $|H|=3$, then $H \in \mathcal{H}$ and we have a different Hamiltonian path in one component, and the third vertex to re-index in another one. If $|H|=2$ then $G[S]=3K_2$ and we are allowed to use the index 3 by the exception in (B4).

We have argued that (B1)-(B4) and (B7) are satisfied during the construction. For (B5), we used both indices $0$ and $f$ only if $G[S] = (f+1)K_1$, and we could pick any two vertices to have $i(s)=0$ and $i(t)=f$. For (B6), note that the number of pairs of indices in $\{k+1, k+2, \dots, f-k-1\}$ is equal to $\binom{f-2k-1}2$ which is more than $f-1 - 2(k+1)$ for each $k=1,2,\ldots,\lceil\frac{f}2\rceil-2$. 

\end{proof}

\vspace{3mm}
\textit{Phase III: Selecting forward orbits}


We first need the following distinction. Call a vertex $x$ \emph{pinned} if it is a neighbour of some vertex $y \in S'<S$ and \emph{unpinned} otherwise (in particular, every vertex of $S$ is pinned). Notice that each pinned vertex has at least one coloured back edge. 


Let $A'$ be the set of these automorphisms in $A$ which fix $S$. 
Fix a vertex $s\in S$ and let $X_1(s),\dots,X_l(s)$, $2 \le |X_1(s)|\le \ldots \le |X_l(s)|$, be all non-trivial orbits with respect to $A'$ of the forward neighbours of $s$. We shall drop the subscript and the argument $s$ if it does not cause confusion or if the choice of $s$ is irrelevant. In particular, the value of $l$ and the cardinalities of these orbits are universal among all the elements of $S$ because $S$ is an orbit with respect to $A$. Define $\beta(X)$ as the number of neighbours in $S$ of any element of a forward orbit $X$ and $p(X)$ as the number of edges between any element of $X$ and $\bigcup \{S'\colon S' <S\}$ (by (A1), all these edges are already coloured). Note that these parameters are well-defined because $X$ is an orbit. Moreover, we shall consider only the orbits $X$ whose cardinalities have to be decreased in order to satisfy (A3), i.e. the ones with $|X| \ge \Delta - p(X) - \beta(X)+1$.

We would like to choose the distribution of the colours on the forward edges in such a way as to split as many orbits as possible. By splitting, we mean that some two vertices in such an orbit get incident edges of different colours to a common neighbour in $S$. Any given vertex $s\in S$ has however a limited splitting ability, as each split requires using two different colours. The number of orbits $s\in S$ can split will be at most $v(s) : = \min\{i(s), f-i(s)\} \le \lfloor \frac{f}2 \rfloor$, and we will refer to this number as the \emph{valency} of $s$. This limitation is necessary so that we keep the vertices of $S$ distinguished.

We would like to satisfy the condition (A3), which requires that after the colouring, each forward orbit $X$ has at most $\max\{1, \Delta-p(X)-\beta(X)\}$ elements. This will be achieved if $X$ is split at least $|X|-\Delta+\beta(X)$ times (if $\Delta-\beta(X)$ is less than one, then the redundant splits will be fictitious -- some two elements of $X$ will get different colours from a single $s\in S$ but it will not change any orbit). Before we describe how to accomplish it, we need some additional observations: 

\begin{itemize}
    \item By \ref{itm:b5}, if $S=S_1$ and $X \subseteq N(z)\cap N(t)$, where $i(z)=f$, $i(t)=0$, then $|X|\leq \Delta-2$, if $S_1$ has two vertices with distinct sets of forward neighbours, or $G = K_{\Delta,\Delta}$ (which satisfies the claim as a Hamiltonian graph) otherwise. 

    \item By the inequality $f\ge |X_1|+\ldots+|X_l|$, if there is an orbit of size at least four (or two orbits of size at least $3$), then $l\le \lfloor\frac{f}2\rfloor -1$. If all orbits but one are of size $2$, the remaining one has size $3$ and $p(X)=0$, then $\beta(X) \ge \Delta-2$ and by counting edges, we obtain that the total number of orbits $X$ is at most $\frac{|S|\cdot f}{2(\Delta-2)} \le f$. 

    \item If $X$ is pinned, then $S\neq S_1$ and for any fixed $k\leq l$, the union $\bigcup_{s\in S} X_k(s)$ must be an orbit, or a union of orbits, with respect to $A$, as otherwise $S$ would not be an orbit with respect to $A$. Therefore, by the procedure assumptions, $|\bigcup_{s\in S} X_k(s)| \le \Delta-p(X_k)$. Note that if $\beta(X_k) < |S|$, then $|X_k| \le \Delta-p(X_k) - |X_k|$, i.e. we can treat $X_k$ as already split at least $|X_k|$ times. 
\end{itemize}

\begin{clm} \label{clm: phase3} Let $i$ be a mapping guaranteed by Claim \ref{clm: i(x)}. It is possible to choose for each vertex $s\in S$ at least $v(s)$ orbits of the forward neighbours of $s$ (with respect to $A'$), so that if every vertex $s$ splits all the chosen orbits, then each forward orbit of $S$ satisfies the condition (A3).
\end{clm}
\begin{proof} 

We consider the vertices of $S$ in any order and for each $s\in S$ such that $v(s)\ge 1$ we decide that we shall split all the orbits of size at least $\lceil\frac{f}2\rceil - v(s) + 2$. It is possible, as the number of such orbits is at most $\frac{f}{\lceil\frac{f}2\rceil - v(s) + 2}$, which proves to be at most $v(s)$ for each $1\le v(s) \le \lfloor\frac{f}2\rfloor$. If $s$ has some spare valency after that, then we split other orbits, starting with the ones of size $2$. In this way, each orbit $X$ is split by all its neighbours of valency at least $\lceil\frac{f}2\rceil - |X| + 2$.

In the following case analysis, we shall check that each orbit has been split enough times.

\textbf{Case 1.} $l=1$. If $\beta(X)=1$ and $p(X)=0$, then any assignment of colours will satisfy (A3). If $\beta(X)=1$ and $X$ is pinned, then either $\beta(X)<|S|$ and $X$ does not require splitting, or $\beta(X)=|S|=1$, and the only element of $S$ has $v(s) = \lfloor\frac{f}2\rfloor \ge \lfloor\frac{|X|}2\rfloor\ge 1$. Assume now that $\beta(X)\ge 2$. By the previous considerations about two vertices with $v(s)=0$ sharing an orbit, we conclude that if there are two vertices not splitting $X$, then $S=S_1$ and $|X|\le \Delta-2$, so $\beta(X)-2$ vertices do split $X$. Otherwise, there are at least $\beta(X)-1$ vertices that split $X$. In both cases, the number of splits is at least $\beta(X)-\Delta+|X|$.

\textbf{Case 2.} $l\ge 2$ and $f\le 4$. Therefore, $f=4$ and there are two orbits $X_1$, $X_2$, each of size two. If $p(X)=0$, then $S\subset N(X)$ and $X$ is split by a vertex with $v(s)=2$. If $X$ is pinned, then either $\beta(X)=|S|$ and the same argument applies, or $\beta(X)<|S|$ and $X$ has a head start of $|X|=2$ free splits which are sufficient with the help of one vertex with $v(s)=2$ or two vertices with $v(s)=1$. 

\textbf{Case 3.} $l \ge 2$, $f\ge 5$ and $|X|\ge 4$. Since $\lceil\frac{f}2\rceil - |X| + 2 \le \lceil\frac{f}2\rceil -2$, then we can use (B6) and argue that $X$ has been split at least $\beta(X) - 2\cdot(\lceil\frac{f}2\rceil - |X| + 2) \ge \beta(X)-f+2|X|-4$ times, which is at least $\beta(X)-\Delta+|X|$ (since $f\le \Delta-1$). 

\textbf{Case 4.} $l \ge 2$, $f\ge 5$,  $|X| \le 3$ and $p(X)=0$. If $\Delta-\beta(X)\le |X|$, then $X$ already satisfies (A3). Therefore, we assume that $\beta(X)\ge \Delta-|X|+1$. This means that $S\subseteq N(X)$, unless $|S|\ge \Delta-1$. In the former case, $X$ is split by a vertex $s_2$ with $v(s)=\lfloor \frac{f}2 \rfloor$. In the later case, consider a vertex $s_3$ with $v(s)=\lfloor \frac{f}2 \rfloor - 1$. By the previous considerations, either:
\begin{itemize}
    \item there is an orbit of size at least four or two orbits of size at least $3$, in which case both $s_2$ and $s_3$ split all their neighbouring orbits (and one of them must be a neighbour of $X$), or
    \item all orbits but one are of size $2$ and the remaining one is of size at most $3$, then there are at most $f$ orbits and $\sum_{s\in S} v(s) \ge 0+0+1+1+\ldots+\lfloor \frac{f}2 \rfloor \ge \lfloor \frac{f}2 \rfloor^2 \ge f+1$, so all the orbits are split once and the one of size $3$ twice.
\end{itemize} 

\textbf{Case 5.} $l \ge 2$, $f\ge 5$,  $|X| \le 3$ and $p(X)>0$. Like in Case 4, $\beta(X)\ge \Delta-|X|+1$, but this time $\beta(X)$ must divide $|S|\le \Delta-1$. For $\Delta\ge 4$ this is only possible if $\beta(X)=|S|$ and then the vertices $s_2$ and $s_3$ split $X$.
\end{proof}

\vspace{\baselineskip}
\textit{Phase IV: Colouring the forward edges}

Let $s\in S$. We already know $i(s)$ and the forward orbits of $s$ (at most $v(s)$ of them) which will be split. In this phase, we shall choose the colours and do the actual splitting. 

We shall choose an \emph{admissible subset} of forward edges of $s$, which consists of $i(s)$ forward edges of $s$ such that for each split forward orbit $X(s)$ we choose at least one, but not all the edges between $s$ and $X(s)$. Such a choice is possible, since the number of these orbits is at most $v(s)$. Recall that we have enumerated the set of colours $C = \{c_0, c_1, \ldots\}$ and that the starting palette is the lowest possible in the lexicographic order. Now, we shall put the lowest possible palette (in the lexicographic order) on the chosen $i(x)$ forward edges. Namely, we put the colour $c_0$ on as many chosen edges as possible, if there is any uncoloured chosen edge left, then we put $c_1$ on as many remaining edges as possible, etc. Finally, over all the admissible subsets, we choose the one, which produces the lowest possible palette (in the lexicographic order), i.e. results in the most edges coloured $c_0$, then $c_1$, etc.

Now, on each other forward edge incident with $s$, we choose the colour with the larger index. 

We now argue that this colouring splits each orbit, which was chosen to be split in Phase III. Let $X(s)$ be a forward orbit of $s$, and assume that all the edges between $s$ and $X(S)$ received the same colour $c_i$. There is an edge $e$ which was not in the chosen admissible subset and the list of this edge must contain, apart from $c_i$, a colour $c_j$ with $j<i$. But that means that interchanging this edge with any other in the chosen admissible subset would result in another admissible subset, which would produce a lexicographically lower palette. This is a contradiction with the choice of the subset.

Unfortunately, at this point, we do not have guaranteed that the vertices with different value of $i$ have distinct forward palettes. We shall slightly modify the current colouring to distinguish these vertices.

Let $s_1, s_2\in S$ be two vertices with $i(s_1)<i(s_2)$ which have the same forward palettes, and assume that there is an automorphism $\varphi$ mapping $s_2$ into $s_1$ (if there is none, we do not have to modify anything). Let $C_1$ be the set of colours on the chosen edges for $s_1$. Let $C_2$ be the set of colours on the images through $\varphi$ of the edges chosen for $s_2$, minus the ones chosen for $s_1$. Finally, let $C_3$ be the set of all other colours of the edges incident to $s_1$. Note that by the choice of the admissible subset, these three sets are non-empty and pairwise disjoint. Moreover, any colour in $C_1$ has a lower index in our enumeration of colours than any colour in $C_2$, and any colour in $C_2$ has a lower index than any colour in $C_3$.

If there is an edge adjacent to $s_1$ coloured with a colour in $C_2$ which has the larger (in our enumeration) colour chosen from its list, then we recolour one such edge with the other colour from its list. Otherwise, we take an edge adjacent to $s_1$ with a colour in $C_2$ and recolour it with the other (larger) colour from its list. Clearly, after this operation, $s_1$ and $s_2$ receive distinct palettes.

If there is any other such pair of vertices, we repeat these steps. We now argue that after a finite number of such operations, we are left with no pair of vertices with identical forward palettes. In the first variant of the operation, we clearly decrease (in the lexicographic order) a palette of some vertex $s_1$. In the second variant, for some vertex $s_1$ and for each possible choice of $s_2$ with $i(s_1)<i(s_2)$, we increase or leave the same (in the lexicographic order) the subpalette of $s_1$ containing the colours from $C_1 \cup C_2$ (and for at least one $s_2$ we increase it). Therefore, after each operation, either the second parameter must increase, or it stays the same while the first parameter decreases. Since both parameters are bounded, this process cannot continue forever.

We are left to argue that we still split the orbits we need. Consider a single operation described in the previous paragraphs. Note that it may have caused the vertex $s_1$ to split one less orbit than before the operation. However, since $s_1$ can be mapped to $s_2$ by an automorphism, this orbit must have been split also by $s_2$. And for $s_2$, it was an additional unplanned split, as it involved colours outside $i(s_2)$. Therefore, each such operation indicates that we initially had a spare split in the required location.

This is a convenient moment to ensure that until this point  the conditions (A1)--(A3) are preserved:

\begin{enumerate}
    \item[(A1)] For each $s\in S$, we coloured all the edges incident to $s$ when we considered that vertex. No other edges have been coloured.
    \item[(A2)] If $S\neq S_1$ and $f=0$, then $S$ must be a singleton -- otherwise by \ref{itm:gwiazdka} for $S$ (which held before the execution of the procedure) $|S|\leq \Delta-b = h < |S|$. If $S\neq S_1$ and $f\ge 1$, then we have broken all the automorphisms of $G[S]$, except possibly the ones that mirror a component of $G[S]$ or the ones that swap the components of $G[S]$. These exceptions are broken by the choice of $i(s)$, $s\in S$, as long as the vertices with different values of $i(s)$ receive different forward palettes, which we already demonstrated above.
    If $S=S_1$, then $f\ge 1$ unless $G=K_n$ and the same argument applies. 
    \item[(A3)] We have argued that this assumption holds during Phase III in the case analysis part. 
\end{enumerate}


\vspace{3mm}
\textit{Phase V: Last recolouring}

Finally, we perform one last modification to ensure that we do not create another vertex (besides $r$ and perhaps the neighbour of $r$) with the starting palette. Assume that, after colouring the edges incident to $S$, there exists a vertex $x$ which satisfies the following conditions:
\begin{itemize}
    \item the vertex $x$ is a forward neighbour of some vertex of $S$,
    \item $x$ has no forward neighbours, 
    \item $x$ is not the root vertex $r$ nor any of the neighbours of $r$, 
    \item in this step we coloured the last back edge incident to $x$, and
    \item the palette of $x$ is the starting palette
\end{itemize}
 Surely, if no such vertex exists, we can continue the algorithm with no risk of violating condition (A4). If, however, such a vertex $x$ exists, let $\hat{X}$ be the orbit of $x$ with respect to $A$ (that is, the orbit just before Phase I of the current iteration of the procedure). 
 
First, we argue that such $x$ must be pinned. 
If it is not, then all the $\Delta$ edges incident to $x$ were coloured in the current step, so $|S|=\Delta$ and $S=S_1$, and one vertex $s\in S$ must have all the forward edges coloured with the larger colours. Since the starting palette was the minimal (in the lexicographic order) among all the vertices, this means that $x$ cannot have it.



Let $\hat{X}_r$ be the set of vertices in $\hat{X}$ which have the same palette as $r$. By minimality of the starting palette, each vertex $x\in \hat{X}_r$ has only the smaller colours chosen from the list of its incident edges. For each such $x$, we choose an incident edge with the largest colour and recolour that edge. This ensures that each such $x$ has a different palette than $r$. However, we might have disrupted the distinction of the elements of $\hat{X}$, as after the change, some automorphism might map a vertex $x\in \hat{X}_r$ to another vertex in $\hat{X}$. Note that the subpalette of such a vertex induced by the edges between $S$ and $\hat{X}$ must be specific, as only one edge may have the larger colour chosen from its list. For brevity, let us call each such edge \emph{big}, and the other edges \emph{small}. For each vertex $x\in \hat{X}$, let $\gamma(x)$ be the number of big edges incident to $x$ in the current colouring. We shall now iteratively (lexicographically) increase the values of $\gamma(x)$ of some vertices that may cause the problem. 

For this, we shall define a sequence $(x_i\colon j\geq 0)$ of distinct elements of $\hat{X}$, and a sequence of their partners $(s_i\colon j\geq 0)$ in $S$. The vertices $x_i$ will be the problematic vertices that may cause the distinguishing issues, and the vertices $s_i$ will be their neighbours in $S$. 

Let $x_0$ be any vertex in $\hat{X}_r$, and assume that we have already defined $(x_i\colon i<j)$ for some $j\geq 1$. If there exists a vertex $x$ in $\hat{X}\setminus \{x_i: i<j\}$, which has at most one big edge with the other endvertex not in  $(s_i \colon i<j)$, then we define $x_j$ as any such vertex with the least possible $\gamma(x)$, and $s_j$ as that other endvertex (if exists). If there is no such vertex, we end the recolourings. 

If we did not define $s_j$ at this step, then $x_j$ must be fixed, and we proceed with the next $j$. Otherwise, let $Y_j$ be the set of vertices in $\hat{X}\setminus \{x_j\}$ which may be mapped into $x_j$. If $Y_j=\emptyset$, then we proceed with the next $j$. If it is not, then for every $y\in Y_j$ and every $i<j$, the edge $ys_i$ (which must exist, as $Y_j$ is a subset of an orbit $\hat{X}$) is big, and there is exactly one more big edge $ys_y$ incident to $y$. Consider these vertices $y\in Y_j$ one by one. If $s_y=s_j$, then we choose a small edge $e$ incident to $y$ with the largest colour, and we change the colour of that edge. If $s_y \neq s_j$, then we increase the colour of $ys_j$. This is the end of the iteration step.

We claim that after Phase V, all the conditions \ref{itm:a1}--\ref{itm:onlyonered} are finally satisfied.
\begin{enumerate}
    \item[(A1)] This condition was satisfied after Phase IV, and in Phase V we recoloured only the edges which were already coloured.
    \item[(A2)] Note that all the elements of $\hat{X}$ must have been fixed after Phase IV, as they had all back edges coloured. Let $x,x'\in \hat{X}$ and assume that there is an automorphism $\varphi$ mapping $x$ to $x'$. This means that at least one of these vertices, say $x$, had an incident edge recoloured in the current phase. In each step $j$, we recoloured only the edges incident to the vertices in $Y_j$, and these vertices had initially $\gamma(y)=\gamma(x_j)$, which was then increased by one. Therefore, each such vertex was later defined as an $x_{j'}$, and in the step $j'$ we broke the automorphism $\varphi$.

Now, consider $s,s'\in S$ such that some automorphism $\varphi$ maps $s$ to $s'$. Again, one of these vertices, say $s$, had an incident edge, say $sx$, recoloured in the current phase. Note that each such vertex $s$ must have been defined at some point as an $s_j$. But each $s_j$ must be fixed, as the unique neighbour of some fixed $x_{j'}$, $j'<j$ through an edge of the largest colour among the edges incident to $x_{j'}$. Therefore, \ref{itm:a2} is also satisfied. 
    \item[(A3)] This follows from (A2), as we have only recoloured edges between the vertices that were fixed before and after Phase V. 
    \item[(A4)] At the beginning of Phase V we recoloured one edge incident to any vertex with the starting palette. Then, we only increased the colours on the edges, so we did not create any new vertex with the starting palette. 
\end{enumerate}

This is the end of the procedure step. By the connectedness of $G$, and since $\le$ is a well ordering,  each vertex of $G$ is reached at some point, and the condition (A2) guarantees that it is fixed. \hfill\qedsymbol

\section{Regular graphs of infinite degree}

In this section, we study the distinguishing edge colourings from lists of regular connected graphs of infinite degree. It was shown in \cite{BP} and \cite{LPSregular} that if $\kappa$ is either equal to $\aleph_0$ or is a fixed point of the aleph hierarchy, then each $\kappa$-regular connected graph has the distinguishing index at most $2$. We extend these results to list colourings.

\mainxx*

\begin{proof}
Let $(v_i \colon i < \kappa)$ be an enumeration of $V(G)$, and $(c_i \colon k<\beta)$ be an enumeration of the colours appearing on the lists for the edges of $G$ for some cardinal $\beta\leq \kappa$. If  $c$ is a colouring from these lists, then we denote by $p_k(v_i)$ the number of the edges incident to $v_i$ coloured with $c_k$.  We define the \emph{palette} of $v_i$ by $(p_k(v_i)\colon k<\beta)$. We shall colour the edges of $G$ in such a way that each vertex shall be given a unique palette. Hence, each vertex of $V(G)$  shall be fixed by every automorphism which preserves the colouring.


We construct a colouring by an induction on $i$. In the step $i$ we colour all the edges between $v_i$ and $v_k$ for every $k>i$, and only these edges. We shall colour them in such a way that the palette of $v_i$ shall be different than the palette of $v_j$ for every $j<i$. 

Fix $i\geq 0$. 
Let $C(j)$ be the number of occurrences of colour $c_j$ on the list of the uncoloured edges incident to $v_i$. First, assume that the set $\mathcal{P}=\{C(j)\colon j<\beta\}$ is unbounded in $\kappa$. It follows that if $\kappa$ is a fixed point of the aleph hierarchy, then there exists a cardinal $\gamma \in \mathcal{P}$ such that $\gamma \geq \aleph_{\aleph_{i^+}}$. If $\kappa=\aleph_0$, then we take $\gamma \in \mathcal{P}$ such that $\gamma>2i+1$. 
Let $c_r$ be a colour for which $C(r)=\gamma$. The set $\{p_r(v_\delta) \colon \delta<i\}$ has at most $|i|$ elements, and there exist more than $2|i|$ cardinal numbers lesser than $\gamma$. Therefore, there exists a cardinal $i<\alpha \leq  \gamma$ such that $\{p_r(v_\delta)\colon \delta<i \}$ does not contain $\alpha$. We colour exactly $\alpha$ of edges incident to $v_i$ with colour $c_r$, including the already coloured edges, and the remaining edges arbitrarily, without using colour $c_r$. Now, $v_i$ is a unique vertex among $(v_j \colon j\leq i)$ which has exactly $\alpha$ incident edges of colour $c_r$. This holds, no matter how the remaining edges of $G$ are coloured.

Now assume that  $\mathcal{P}$ is not unbounded in $\kappa$. It follows that $\beta=\kappa$, and there exists $\kappa$ colours on the lists of the edges incident to $v_i$. Let $\mathcal{S}=\{j \colon j < \kappa)$ be an enumeration of these colours. We conduct an induction on $\delta<i$. Assume that less than $\kappa$ edges incident to $v_i$ has  already been coloured. In the step $\delta$, we pick a colour $c_r$ which does not yet appear on the edges incident to $v_i$.  If $c_r$ does not belong to the palette of $v_\delta$, then we colour  exactly one edge incident to $v_i$ with $c_r$, and the remaining uncoloured edges incident to $v$ which have $c_r$ on their lists with different colours. If not, then we do not use colour $c_r$ on any edge incident to $v_i$, and we colour each edge incident to $v_i$ with colour $c_r$ in its palette with colour different than $c_r$. It follows that $p_r(i)\neq p_r(\delta)$ no matter how other edges of $G$ shall be coloured. Hence, $v_i$ and $v_j$ have different palettes. After step $\delta$, we coloured less than $\kappa$ edges incident to $v_i$ because $\mathcal{P}$ is not unbounded in $\kappa$. Hence, the assumption for the induction on $\delta$ always holds. After the induction on $\delta$, we colour the remaining uncoloured edges incident to $v_i$ arbitrarily.

\end{proof}

We conjecture that the statement of the theorem above holds for every cardinal number $\kappa$ i.e. for every connected regular graph of infinite degree. We formulate this conjecture in a form similar to Theorem \ref{thm:aleph}.

\begin{conj}
  Let $\kappa$ be an infinite cardinal number, and $G$ be a connected  $\kappa$-regular graph. Then $D'_l(G)\leq 2$.
\end{conj}

\bibliographystyle{abbrv}
\bibliography{lit.bib}

\end{document}